
\input amstex

\input amsppt.sty
\magnification=\magstep1
\hsize15truecm
\vsize24truecm

\font\midbf=cmbx10 scaled1200

\def\QED{\hfill{\bf Q.E.D.}}
\def\c_t{\chi_{\text{top}}}
\def\t{\widetilde}
\def\R{\Cal R}
\def\B{\Cal B}
\def\wt{\widetilde}
\def\la{\leftarrow}
\def\o{\over}

\def\G{\Gamma}

\font\eightrm=cmr10 at 8pt

\font\bigbf=cmbx10 scaled1440
\font\midbf=cmbx10 scaled1200

\title\nofrills $ $ \\ \\ \\
{\bigbf On~the~invariants~of~base~changes \\
 of~pencils~of~curves, II }\endtitle

\leftheadtext\nofrills{S.-L. Tan}
\rightheadtext\nofrills{Invariants of base changes, II}

\author{\hskip0cm \midbf Sheng-Li Tan $^\star$ } \endauthor

\address
 Department of Mathematics,
 East China Normal University,
 Shanghai 200062,
 P.~R.~of China\endaddress

\curraddr
 Max-Planck-Institut f{\"u}r Mathematik,
 Gottfried-Claren-Str. 26,
 53225 Bonn, Germany \endcurraddr

\thanks$^\star$
The author would like to thank the hospitality and financial support of
Max-Planck-Institut f{\"u}r Mathematik in Bonn during this research.
This research is partially supported by the National Natural Science
Foundation of China  and by the
Science Foundation of the University Doctoral Program of CNEC.  \endthanks

\endtopmatter
\hsize15truecm
\vsize24truecm
\baselineskip 14.3truept
\vskip1.2cm

\noindent{\bf Introduction}
\vskip0.6cm
\noindent
    Semistable reductions of pencils of curves have been studied
by many authors in various ways. (cf. [AW, De, DM, Xi3]). In this part
of the series, we shall investigate semistable reductions from the point of
view of numerical invariants. As an application, we obtain two numerical
criterions
for a base change to be stabilizing, and for a fibration to be isotrivial.
We also obtain a canonical class inequality for any fibration.
Some other applications are presented.

   Let $f:S\longrightarrow C$ be a fibration of a smooth complex projective
surface $S$ over a curve $C$, and denote by $g$ the genus of a general fiber of
$f$. We assume that $g>0$ and    $S$ is relatively minimal with respect to $f$,
i.e.,
$S$ has no $(-1)$-curves contained in a fiber of $f$. The basic relative
numerical invariants of $f$ are defined as follows,
   $$
      \eqalign{ \chi_f&=\chi(\Cal O_S)-(g-1)(g(C)-1), \cr
      K_f^2&=K^2_S-8(g-1)(g(C)-1), \cr
      e_f&=\chi_{\text{top}}(S)-4(g-1)(g(C)-1).}
   $$
These invariants are nonnegative integers satisfying the Noether equality
 $12\chi_f=K^2_f+e_f$. We denote by $\omega_{S/C}=\omega_S
\otimes f^*\omega_C$ the relative canonical sheaf of $f$, and $K_{S/C}$
the relative canonical divisor corresponding to $\omega_{S/C}$. Then
$\chi_f=\deg f_*\omega_{S/C}$ and $ K_f^2=K_{S/C}^2$.
If $g>1$ and $f$ is not locally
trivial, then $\chi_f$ and $K_f^2$ are positive (cf. [Ar, Be2, Pa] or
[BPV, Theorem~18.2] and [Xi1, Theorem~2]), in this case, we define the slope of
$f$ as
               $$
                 \lambda_f=K_f^2/\chi_f.
               $$
$e_f=\sum_{F}e_F=\sum_{F}(\chi_{\text{top}}(F)-(2-2g))$ is zero iff
$f$ is smooth.

A fiber of $f$ is called semistable if it consists of simple components meeting
normally. $ f$ is said to be semistable if every fiber of it is semistable.

   Let $\pi: \widetilde C\longrightarrow C$ be a base change of degree $d$.
 Then the pullback fibration
$\widetilde f:\widetilde S\longrightarrow \widetilde C$ of $f$ with respect to
$\pi$
is defined as the relative minimal model of the desingularization of
$S\times_C\widetilde C \longrightarrow \widetilde C$. (cf. Sect.~1.3).
Since $g>0$, the relative minimal model is unique, hence $f$ is determined
uniquely by $f$ and $\pi$. Due to Kodaira's classification of singular
fibers, the semistable reduction of an elliptic fibration is quite clear, so we
always assume that $g\geq 2$.

We  define
  $$
    \chi_\pi=\chi_f-{1\over d}\chi_{\widetilde f}, \hskip0.3cm
    K_\pi^2=K^2_f-{1\over d}K^2_{\widetilde f}, \hskip0.3cm
    e_\pi=e_f-{1\over d}e_{\widetilde f}
  $$
as the basic numerical invariants of $\pi$ with respect to $f$. Obviously,
they are rational numbers satisfying $12\chi_\pi=K^2_\pi+e_\pi$.
Xiao [Xi4] and I [Ta1] proved that these invariants are nonnegative,
and one of them vanishes if and only if $\pi$ is an {\it invariant base
change}.
(See Definition 1.7).

\proclaim{Definition~I} {\rm We shall call $\pi$ a {\it stabilizing}
(resp. {\it trivial}) base change if all of the fibers of $\widetilde f$ (resp.
$f$)
over the branch locus $B_\pi$
of $\pi$  are semistable. We shall also call $\pi$ the {\it semistable
reduction of the fibers over $B_\pi$}. }
\endproclaim

     The well-known semistable reduction theorem says that for any fibration
$f$,
there exists a base change $\pi$ such that $\t f$ is semistable.
In particular,
  let $\pi$ be a base change totally ramified over $F$ (i.e., over $f(F)$) and
some other
semistable fibers, and let $F'$ be the minimal embedded
resolution of $F$. If the degree of $\pi$
is exactly the least common multiple of the multiplicities of the
components in $F'$, then it is well-known
that $\pi$ is stabilizing. We shall call $\pi$ the {\it canonical semistable
reduction } of $F$,
and denote it by $\phi_F$.

\proclaim{Definition~II} {\rm For any fiber $F$ of $f$, we define its basic
invariants
to be the basic invariants of $\phi=\phi_F$, and denote them respectively by
   $$
   c^2_1(F)=K_{\phi}^2, \hskip0.4cm
   c_2(F)=e_{\phi},  \hskip0.4cm
   \chi_{F}=\chi_{\phi}.
   $$}
\endproclaim

We shall show that these invariants are independent of the  choice of
the base changes (Lemma~2.3). They are nonnegative rational numbers
satisfying  the Noether equality
     $$12\chi_F=c_1^2(F)+c_2(F).$$
We can also see that one of them vanishes iff $F$ is semistable (Lemma~1.8).
In fact, these invariants can be computed directly from the embedded resolution
of $F$
(see Theorem~3.1 for the formulas). For simplicity, if $B=F_1+\cdots+ F_s$,
then we define $c_1^2(B)=c_1^2(F_1)+\cdots+ c_1^2(F_s)$. Similarly, we can
define $c_2(B)$ and $\chi_B$.

\proclaim{Definition~III} {\rm A fibration $f: S\longrightarrow C$ is {\it
trivial} if $S$ is isomorphic
 to $F\times C$ over $C$. It is {\it isotrivial} if it becomes trivial after a
finite base change.}
\endproclaim

If $\widetilde f$ is a semistable model of $f$ under a semistable reduction
$\pi$,
then
a natural problem is:

\vskip0.2cm

{\it What is the effect of a non-semistable fiber on the invariants  of $\t
f$~? }
\vskip0.2cm
\noindent
(See [Xi2, Problem 7]). In this paper the effect is completely determined.

In what follows, we denote
by $\B_\pi=f^*(B_\pi)$  the locus of branched fibers, and by
$\R_\pi=\t f^*(R_\pi)$  the pullback fibers of $\Cal B_\pi$, where
$R_\pi$ is the inverse image of $B_\pi$ under $\pi$.

The main results of this paper are the following.

\proclaim{Theorem~A} Let $f: S\longrightarrow C$ be a fibration, and let
$\pi: \t C\longrightarrow C$ be a base change of degree $d$. Then
    $$
      K_\pi^2=c_1^2(\B_\pi)-{1\o d}c_1^2(\R_\pi),\hskip0.3cm
      e_\pi  =c_2(\B_\pi)  -{1\o d}c_2(\R_\pi),\hskip0.3cm
     \chi_\pi=\chi_{\B_\pi}-{1\o d}\chi_{\R_\pi}.
    $$
\endproclaim

\proclaim{Corollary} For any fibration $f:S\longrightarrow C$ and any base
change
 $\pi: \widetilde C\longrightarrow C$, we have

{\rm 1)} $$
      K_\pi^2 \leq c_1^2(\B_\pi), \hskip0.5cm
      e_\pi  \leq c_2(\B_\pi) , \hskip0.5cm
     \chi_\pi \leq \chi_{\B_\pi},
    $$
and one of the equalities holds iff $\pi$ is stabilizing.

{\rm 2)}  $$\sum_{F}c_1^2(F)\leq K_f^2, \hskip0.3cm
      \sum_F \chi_F\leq \chi_f, \hskip0.3cm
      \sum_F c_2(F)\leq e_f,
    $$
where $F$ runs over all of the non-semistable fibers of $f$. Furthermore,
one of the first two equalities holds iff $f$ is isotrivial, and the last
equality holds iff the semistable model of $f$ is smooth.

{\rm 3)} \
 If $f$ is non-isotrivial and $\pi$ is stabilizing, then we have
   $$\lambda_{\t f}={K_f^2-c_1^2(\B_\pi)\o \chi_f-\chi_{\B_\pi}}.$$
Hence the slope of $\t f$ is completely determined by the branched
non-semistable fibers.
\endproclaim

{}From the point of view of the fibration itself, we can define
    $$
       I_K(f)=K^2_f-\sum_{F}c_1^2(F), \hskip0.2cm
       I_\chi(f)=\chi_f-\sum_F\chi_F, \hskip0.2cm
       I_e(f)=e_f-\sum_{F}c_2(F),
    $$
where $F$ runs over all of the singular fibers of $f$. Then we know that
these numbers are nonnegative invariants of $f$, and the first two
invariants measure the isotriviality of $f$. (Note that $K^2_f$ and $\chi_f$
measure the local triviality of $f$). If $f$ is semistable,
then these invariants are nothing but the standard relative invariants
of the fibration. Furthermore, if we think of $\R_\pi$ (resp. $\B_\pi$) as
the set of singular fibers of $\widetilde f$ (resp. $f$), then it is easy
to see that Theorem~A holds too. Hence Theorem~A can be restated as

\proclaim{Theorem~A$'$} For any fibration $f$ and any base change $\pi$ of
degree $d$, we have
    $$
       I_K(\wt f)=dI_K(f),\hskip0.2cm
       I_\chi(\wt f)=dI_\chi(f),\hskip0.2cm
       I_e(\wt f)=dI_e(f).
    $$
\endproclaim

Due to Theorem~A, the study of the invariants of stabilizing base changes
can be reduced to the local study of $c_1^2(F)$ and $c_2(F)$. First of all,
from definition, it is trivial to see that
    $$ c_2(F)\leq e_F\ (=:\chi_{\text{top}}(F)-(2-2g)),$$
with equality iff the semistable model of $F$ is a smooth fiber.
In Sect.~3.3, we obtain

\proclaim{Theorem~B}
    $$ c_1^2(F)\leq 2c_2(F), $$
with equality iff $F=nF_{\text{red}}$ and $F_{\text{red}}$ has at worst
ordinary double points as its singularities. Hence for any stabilizing
base change $\pi$, we have
     $$ K_\pi^2\leq 8\chi_\pi.$$
In particular, if $S$ admits an isotrivial fibration, then
     $$K^2_S\leq 8\chi(\Cal O_S),$$
with equality iff all of the singular fibers are some multiples of smooth
curves.
\endproclaim

We show that $c_1^2(F)$ is in fact bounded by the genus $g$, i.e.,
\proclaim{Theorem~C}
         $$
             c_1^2(F)\leq 4g-4.
         $$
\endproclaim

As an application of this inequality,
we obtain the following {\it canonical class inequality}.

\proclaim{Theorem~D} If $f$ is a non-trivial fibration of genus $g\geq 2$ with
$s$
singular fibers, then
    $$
       K_{S/C}^2\leq (2g-2)(2g(C)-2+3s),
     $$
and if the equality holds, then $f$ is smooth, i.e., $s=0$.
\endproclaim

Note that other canonical class inequalities are already known for non-trivial
semistable fibrations:
    $$\eqalign {&K_{S/C}^2\leq (2g-2)(2g(C)-2+s);\cr
               &K_{S/C}^2< 4g(g-1)(2g(C)-2+s);\cr
               &K_{S/C}^2\leq 8(g-1)^2(2g(C)-2+s).}
    $$
These inequalities are due respectively to Vojta [Vo], Szpiro [Sz], Esnault and
Viehweg
[EV]. Note that in [Ta2], we have shown that the equality in Vojta's inequality
implies the smoothness of $f$.

As another application, we find some new phenomena for fibrations. (Sect.~4.1).
For example, from the corollary above, we can see that every
semistable model $\t f$ of
$f$ has the same slope
    $$
         \lambda_{\wt f}={I_K(f)/ I_\chi(f)}.
    $$
{}From Theorem~B we know that if $\lambda_f>8$, then any
non-trivial stabilizing base change
$\pi$ makes the slope increase.
We have also found some relationships
between non-semistable fibers and the slope of a fibration.

Finally, in Sect.~4.3, we consider the computation of the Horikawa
number of a genus 3 non-semistable fiber $F$ through semistable reductions.
We reduce it to the computation  for its semistable models $\t F$.

\proclaim{Notations} {\rm
If $D$ is a local curve and $p\in D$, then we denote by $\nu_p$ the
multiplicity of $D$ at $p$, and denote respectively by $\mu_p$,
$\delta_p$, $k_p$ the Milnor number, geometric genus and the number
of local branches of $(D_{\text{red}}, p)$. Hence $\mu_p=2\delta_p-k_p+1$.
If $F$ is a curve on a smooth surface, then we denote by $\mu_F$
the total Milnor number of the singularities of $F$.

If $a, b$ are two natural numbers, then we denote by $(a, b)$ the greatest
common divisor of $a$ and $b$, and let $[a, b]={(a, b)^2\o ab}$. $[x]$ is
the greatest integer $\leq x$
}
\endproclaim

\vskip0.8cm
\noindent{\bf 1 \ Preliminaries and technical lemmas }
\vskip0.6cm
\noindent
{\it 1.1 \ Embedded resolution of curve singularities}
\vskip0.6cm
\noindent
Let $(B,p)\subset \Bbb C^2$ be a local curve (not necessarily reduced) in a
neighborhood $U_0$ of $p=(0,0)$. Assume that $(B_{\text{red}},p)$ is a singular
point, we say also that $p$ is a singular point of $B$.

\proclaim{Definition 1.1} {\rm  The embedded resolution of a curve singularity
$(B, p)=(B_0, p_0)$ is a sequence
   $$
      (U_0, B_0)\overset{\sigma_1}\to\la
      (U_1, B_1)\overset{\sigma_2}\to\la
      \cdots
       \overset{\sigma_r}\to\la
       (U_r, B_r)
   $$
satisfying the following conditions.

    (1) $\sigma_i$ is the blowing-up of $U_{i-1}$ at a singular point
$p_{i-1}\in B_{i-1}$ with $\mu_{p_{i-1}}>1$.

    (2) $B_{r, \text{red}}$ has at worst ordinary double points
             as its singularities.

    (3) $B_i$ is the total transformation of $B_{i-1}$.
}
\endproclaim

It is well-known that embedded resolution exists and is unique for any
curve singularity $(B,p)\subset \Bbb C^2$.

We denote by $m_i$ the the multiplicity of $(B_{i,{\text{red}}},p_i)$. Let
    $$
       \alpha_p=\sum_{i=0}^{r-1}(m_i-2)^2.                            \eqno(1)
    $$
If $q\in B_r$ is a double point, and $a_q, b_q$ are the multiplicities
of the two components of $(B_r, q)$, then we let
    $$
       \beta_p=\sum_{q\in B_r}[a_q, b_q].                             \eqno(2)
    $$

\proclaim{ Lemma~1.2}
        $$
   \eqalignno{
     \mu_p &=\sum_{i=0}^{r-1}(m_i-1)(m_i-2)+k_p-1,                   &(3)
     \cr
  \delta_p&={1\o2}\sum_{i=0}^{r-1}(m_i-1)(m_i-2)+k_p-1   .           &(4) }
       $$
\endproclaim
\demo{Proof}  In the embedded resolution, we let $E_1\cap (B_1-E_1)=p_1,
\cdots, p_s$. Then by [Ta1, Lemma~1.3] we have
    $$
      \mu_p=(m_p-1)(m_p-2)-1+\sum_{i=1}^s\mu_{p_i}.                  \eqno(5)
    $$
On the other hand, it is obvious that
    $$
      k_p=\sum_{i=1}^s(k_{p_i}-1),                                 \eqno(6)
    $$
hence (3) can be obtained easily by using induction on $r$, and (4)
follows from (3) and $\mu_p=2\delta_p-(k_p-1)$.
\QED
\enddemo

\proclaim{Lemma~1.3} For any singular point $(B, p)$, we have
    $$
    \alpha_p+\beta_p \leq \mu_p.                                    \eqno(7)
    $$
\endproclaim
\demo{Proof} First we prove (7) for the case $m_p=2$, i.e., $(B_{\text{red}}
, p)$ is a double point. Assume that $(B,p)$ is defined by $f(x,y)=0$ at $0$.

    If $f=x^a(x+y^k)^b$ and $k=1$, then $\alpha_p=0$, $\mu_p=1$ and
$\beta_p=[a, b]$,
(7) is obvious. If $k>1$, then by the computation of  the embedded resolution,
we have
    $$
       \alpha_p=k-1, \ \mu_p=2k-1, \ \beta_p=1-{1\o k}+[a, k(a+b)]+[b,
k(a+b)]\leq 1,
    $$
hence (7) holds strictly.

    If $f=(x^2+y^{2k+1})^n$, then
    $$
       \alpha_p=k, \ \mu_p=2k, \ \beta_p={3\o 2}(1-{1\o 2k+1}),
    $$
thus we can see  that $\alpha_p+\beta_p\leq \mu_p$.

Now we assume that $m_p\geq 3$. In this case, we shall prove (7)
 by using induction on $\mu_p$. From (5) we know
 $\mu_{p_i}<\mu_p$, by induction hypothesis, we have $\alpha_{p_i}
+\beta_{p_i}\leq \mu_{p_i}$. On the other hand, we know
    $$
  \beta_p=\sum_{i=1}^s\beta_{p_i}, \hskip0.3cm \alpha_p=(m_p-2)^2+\sum_{i=1}^s
\alpha_{p_i},
    $$
from (5), (7) follows immediately.
\QED
\enddemo

\vskip0.8cm
\noindent
{\it 1.2 \ On the  resolution of the singularity of $z^d=f(x,y)$}
\vskip0.6cm
\noindent
Now we assume that $(B,p)$ is defined by $f(x,y)=0$ at $p=(0,0)$. Let $\Sigma
\subset \Bbb C^3$ be a local surface defined by $z^d=f(x,y)$, and let
$V_0$ be the normalization of $\Sigma$. Then, $V_0$ is a $d$-cyclic cover
$\pi_0:
V_0\longrightarrow U_0$, the singular points of $V_0$ (lying over $p$) can be
resolved by
the embedded resolution of $(B, p)$, it goes as follows.

Let $V_r$ be the normalization of $U_r\times_{U_0}V_0$, and let
$\eta : M\longrightarrow V_r$ be the minimal resolution of the singularities of
$V_r$.

    $$
         \CD
          V_0@<\tau<<V_r@<\eta<<M\\
          @V\pi_0VV@V\pi_rVV@VV{\pi_r\eta}V \\
          U_0@<<\sigma<U_r@=U_r
         \endCD
    $$
Then $\pi_r$ is a cyclic covering branched along $B_r$. If near $q\in
B_r$, $B_r$ is defined by $x^ay^b=0$, then $V_r$ is locally the normalization
of
$z^d=x^ay^b$, which are cyclic quotient singularities, hence can be resolved
by Jung-Hirzebruch method (cf. [BPV, p.83]). Hence
$\phi=\tau\eta:M\longrightarrow V_0$
is the resolution of $V_0$, we shall call $\phi$ the {\it embedded
resolution of $V_0$}.

    Denote by $E_{p}=\sum_{i=1}^sE_i$ the exceptional curves of $\phi$,
and let $K_\phi=\sum_{i=1}^sr_iE_i$ be the rational canonical divisor
of $E_{p}$, which is determined uniquely by the adjunction formula
$K_\phi E_i+E_i^2=2p_a(E_i)-2$. Then $K_{\phi}^2$ is an invariant
of the resolution $\phi$. If $\phi$ is minimal, then $K^2_\phi=K^2_{p}\leq 0$
is an invariant of the singularities of $V_0$, which is independent
of the resolution. $K_{p}^2=0$ iff $V_0$ has at worst rational
double points as its singularities. At the end of this paper, we shall present
Jung-Hirzebruch's resolution and the computation of $K^2$ for the singularity
defined by $z^d=x^ay^b$.
The following lemma can also be obtained from this resolution. (cf. Sect.~5,
or [Xi3] and [BPV, p.83]).

\proclaim{Lemma~1.4} If $(B,p)$ is defined by $x^ay^b=0$, and $d$ is divided by
$a$ and $b$, then $E_p$ consists of $d_p=(a,b)$ disjoint curves of type $A_n$,
where
  $$
     n=[a, b]{d\over d_p}-1 .                       \eqno(8)
  $$
\endproclaim

\proclaim{Lemma~1.5} Assume that $d$ is divided by all of the multiplicities of
the components in the embedded resolution $B_r$. Then
   $$
    -{1\o d}K^2_\phi=\alpha_p.                                    \eqno(9)
   $$
\endproclaim

The proof of this lemma will be given in Sect.~5.

Now we recall the normalization of $\Sigma$. (cf. [Ta1, Lemma~2.1]).

\proclaim{Lemma~1.6} For any point $p\in B$, $\pi_0^{-1}(p)$ consists of
$d_{p}=\gcd(d, n_1, \cdots, n_s)$ points if there are exactly $s$ components
$\G_1, \cdots, \G_s$ passing through $p$, where $n_i$ is the multiplicity of
$\G_i$ in $B$.
\endproclaim

\vskip0.8cm
\noindent
{\it 1.3 \ The construction of base changes}
\vskip0.6cm
\noindent
In this section, we recall the construction of the pullback fibration $\t f$
of $f:S\longrightarrow C$ under a base change.

   Let $\pi: \widetilde C\longrightarrow C$ be a base change of degree $d$.
Then the pullback fibration
$\widetilde f:\widetilde S\longrightarrow \widetilde C$ of $f$ with respect to
$\pi$
is defined as the relative minimal model of the desingularization of
$S\times_C\widetilde C \longrightarrow \widetilde C$.
In fact, the pullback fibration $\widetilde f:\widetilde S
\longrightarrow \widetilde C$ can be constructed as follows.

    Let $\rho_1:S_1\longrightarrow S\times_C\widetilde C $ be the normalization
of $S\times_C\widetilde C$, let $\rho_2:S_2\longrightarrow S_1$ be the
minimal desingularization of $S_1$. Then we have a fibration
$f_2:S_2\longrightarrow \widetilde C$. Let $\widetilde \rho:
S_2\longrightarrow  \widetilde S$ be the contraction of $(-1)$-curves such that
$\widetilde f:\widetilde S\longrightarrow \widetilde C$ is a
relative minimal model.
Since we have assumed that $g>1$, $\widetilde \rho$ is unique.
Hence $\t f:\t S\longrightarrow \t C$ is determined uniquely by $f$ and $\pi$.

   $$\CD
     \widetilde S@<\widetilde\rho<< S_2 @>\rho_2>> S_1 @>\rho_1>>
                              S\times_C\widetilde C @>\Pi'>> S \\
      @VV\widetilde fV @VVf_2V @VVf_1V @VVV @VVfV \\
       \widetilde C @=\widetilde C@=\widetilde C@=\widetilde C@>>\pi> C
      \endCD
    $$
Let $\Pi_2=\Pi'\circ\rho_1\circ\rho_2:S_2\longrightarrow S.$

\proclaim{\noindent Definition 1.7} {\rm If $\pi:\widetilde C\longrightarrow C$
is a base change satisfying
       $$\widetilde \rho^*K_{\widetilde S/\widetilde C}\equiv \Pi_2^*K_{S/C},
$$
then we shall call it an {\it invariant  base  change}}.
\endproclaim

In fact, we have shown that if $g\geq 2$, then $\pi$ is invariant iff the
fibers $F$ over the branch
locus are reduced and $F$ has at worst $d_F$-simple singularities,
where $d_F$ is the greatest ramification index of $\pi$ over $f(F)$.
(cf. [Ta1, Lemma~2.2 and 2.3]).
A $d$-{\it{simple singularity}} is a simple curve singularity $f(x, y)=0$ such
that $z^d=f(x, y)$ is a simple surface singularity. Hence $2$-simple is $ADE$,
3-simple is $A_1$, $\cdots$, $A_4$, 4 and 5-simple are
$A_1, A_2$, $d$-simple is $A_1$ if $d>5$.

    In particular, we can see that the canonical semistable reduction $\phi_F$
of $F$
is invariant if and only if $F$ is semistable. On the other hand, In
[Xi4, Ta1] we have proved that the basic invariants of a base change are
nonnegative, and one of them vanishes iff the base change is invariant. Hence
we have

\proclaim{Lemma~1.8} If $g\geq 2$, then the invariants $c_1^2(F)$, $c_2(F)$ and
$\chi_F$ are nonnegative, and one of them vanishes if and only if $F$ is
semistable.
\endproclaim

Let $F$ be a singular fiber. We always denote by $F'$ the
embedded resolution of $F$,
and denote by $M_F$ the least common multiple of the multiplicities
of the components in $F'$.

\vskip0.8cm
\noindent
{\bf 2 \ On the invariants of a base change }
\vskip0.6cm
\noindent
{\it 2.1 \ Local computations of $K^2_\pi$ }
\vskip0.6cm
\noindent
     In this section, we first consider the computation of the invariant
$K^2_\pi$
for a base change $\pi: \widetilde C\longrightarrow C$. Without loss of
generality, we
assume that $\pi$ is totally ramified over $p_1, \cdots, p_s$. Let
$\rho_2$ be the embedded resolution of singularities, let $F_1, \cdots, F_s$
be the fibers of $f$ corresponding to $p_1, \cdots, p_s$, and let $\Cal B_\pi
=\sum_{i=1}^sF_i=\sum_\Gamma n_\Gamma\Gamma$. From Lemma~1.6, it is easy to see
that
    $$
       K_{S_2}\equiv\Pi_2^*\left(K_S+\sum_{\G\subset \Cal B_\pi}
       \left(1-{(d, n_\G)\over d}\right)\G\right)
              +K_{\rho_2}.                                \eqno(10)
    $$
where $K_{\rho_2}$ is the rational canonical divisor of the exceptional
set of $\rho_2$. On the other hand, we have
    $$
       K_{\widetilde C}=\pi^*\left(K_C+\sum_{i=1}^s\left(1-{1\over d}\right)
                                                  p_i\right).    \eqno(11)
    $$
Note that $f_2^*\pi^*=\Pi_2^*f^*$, hence from (10) and (11) we can obtain

    $$
      K_{S_2/\widetilde C}=\Pi_2^*\left(K_{S/C}-\sum_{i=1}^sH_{F_i}\right)
                  +K_{\rho_2},
    $$
where $H_i=\sum_{\G\subset F_i}h_\G\G$, $h_\G=n_\G-1-{1\over d}(n_\G-(d,
n_\G)).$
Hence
    $$
      dK_f^2-K_{f_2}^2=d\sum_{i=1}^s(2H_{F_i}K_S-H_{F_i}^2)-K_{\rho_2}^2.
    $$
If we let $K^2_\pi(f_2) =K_f^2-{1\o d}K^2_{f_2}$, then
    $$
      K_\pi^2=K_\pi^2(f_2)-{1\o d}\#\{\, (-1)\text{-curves contracted by
}\widetilde
            \rho\, \}.
    $$
\proclaim{Proposition~2.1} With the notations above, we have
    $$
       K_\pi^2(f_2)=\sum_{i=1}^s(2H_{F_i}K_S-H_{F_i}^2)-
              \sum_{i=1}^s\sum_{p\in F_{i}}{1\o d}K_p^2,           \eqno(12)
    $$
\endproclaim

In the case when $\pi$ is the base change of $F_1, \cdots, F_s$,  if
$d$ is divided by  $M_{F_i}$, $i=1, \cdots s$,
we can see that $H_{F_i}=F_i-F_{i, \text{red}}$. Note that we have (cf. [Ta1,
(7)])
     $$e_F=\chi_{\text{top}}(F)-(2-2g)=2N_F+\mu_F, $$
where $N_F=g-p_a(F_{\text{red}})
={1\o2}((F-F_{\text{red}})K_S-F_{\text{red}}^2)$ is an invariant of $F$.
{}From Lemma~1.5 we have

\proclaim{Proposition~2.2} If $d$ is divided by $M_{F_i}$ for all $i$, then
    $$
   K_\pi^2(f_2)=\sum_{i=1}^s(4N_{F_i}
+F_{i, \text{red}}^2)+\sum_{i=1}^s\sum_{p\in F_i} \alpha_p.          \eqno(13)
    $$
\endproclaim

Note that the right hand side of (13) is independent of $d$.

\vskip0.8cm
\noindent
{\it 2.2 \ Proof of Theorem~A}
\vskip0.6cm
\noindent
We consider first the composition of base changes.

Let $\pi_1:C_1\longrightarrow C$ and $\pi_2: \t C\longrightarrow C_1$ be two
base changes, let
$f_1$ be the pullback fibration of $f$ under $\pi_1$, and let $f_2$ be that
of $f_1$ under $\pi_2$. By the universal property of fiber product and
the uniqueness of the relative canonical model (when $g>0$), we know
$f_2$ is nothing but the pullback fibration $\t f$ of $f$ under
$\pi=\pi_1\circ\pi_2$. Hence we have the {\it basic equalities}:
    $$\eqalign{ K_\pi^2&= K_{\pi_1}^2+{1\o \deg\pi_1}K_{\pi_2}^2, \cr
                e_\pi  &= e_{\pi_1}+{1\o \deg \pi_1}e_{\pi_2}, \cr
               \chi_\pi&= \chi_{\pi_1}+{1\o \deg\pi_1}\chi_{\pi_2}.}
                                                                      \eqno(14)
    $$

\proclaim{Lemma~2.3} Let $f: S\longrightarrow C$ be a fibration, and let $F_1,
\cdots, F_s$
be fibers of $f$. Considering all of the semistable reductions $\pi$ of $F_1,
\cdots,
F_s$, we have that $K_\pi^2$, $e_\pi$ and $\chi_\pi$ are independent of $\pi$.
\endproclaim
\demo{Proof} Let $\pi_1:C_1\longrightarrow C$ and $\pi_2:C_2\longrightarrow C$
be two semistable reductions
of $F_1, \cdots, F_s$,  let $\deg\pi_i=d_i, $ $i=1, 2$, and  let
$f_i$ be the pullback fibration of $f$ under $\pi_i$. We shall prove that
        $$K_{\pi_1}^2=K_{\pi_2}^2.$$

     For this, we consider the pullback of $\pi_1$ and $\pi_2$,
    $$
         \pi=\pi_1\times_C\pi_2: \t C=C_1\times_CC_2\longrightarrow C.
    $$
Note that if necessary, we can choose $\widetilde C$ to be the normalization
of a component of $C_1\times_CC_2$. Let $p_i:\t C\longrightarrow C_i$ be the
$i$-th projection,
 it is obvious that
$\deg p_1=d_2$, $\deg p_2=d_1$. Then we have
$\pi=p_1\circ\pi_1=p_2\circ\pi_2$ (composition of base changes). $\pi_1$
and $\pi_2$ are semistable reductions, so the fibers of $f_i$ over $F_1,
\cdots, F_s$ are semistable,
and thus $p_i$ is an invariant base change.
It implies that $K_{p_i}^2=0$ for $i=1, 2$.
Then by using the basic equalities (14), we have
    $$
        K_\pi^2=K_{\pi_1}^2=K_{\pi_2}^2.
    $$

The proofs for $\chi_\pi$ and $e_\pi$ are the same as above.
\QED
\enddemo

\proclaim{Lemma~2.4} In the situation of Lemma~2.3, we have
   $$K_\pi^2=\sum_{i=1}^sc_1^2(F_i), \hskip0.2cm
     e_\pi  =\sum_{i=1}^sc_2(F_i), \hskip0.2cm
     \chi_\pi=\sum_{i=1}^s\chi_{F_i}.
   $$
\endproclaim
\demo{Proof} By Lemma~2.3, we can assume that $\pi$ is the pullback of
the canonical semistable reductions $\pi_i=\phi_{F_i}:C_i\longrightarrow C$,
$i=1, \cdots, s$.
We can assume that $\pi_i$ is unramified over the fibers $F_j$ for $j\neq i$.
Without loss of generality, we assume also that $s=2$. As in the proof of
Lemma~2.3,
we have
     $$K_\pi^2=K_{\pi_1}^2+{1\o d_1}K_{p_1}^2.$$
Since $p_1$ is a totally ramified semistable reduction, $K_{p_1}^2$ can be
computed locally from the branched  non-semistable fibers, which are the
pullback of $F_2$ under $\pi_1$. Hence we know
     $$K_{p_1}^2=d_1K_{\pi_2}^2. $$
By definition, $K_{\pi_i}^2=c_1^2(F_i)$. Hence
we have obtained the desired equality.

Note that the local property used above holds for $e_{\pi}$ and
$\chi_\pi$.
\QED
\enddemo

\demo{Proof of Theorem~A}
Let $\hat \pi:\hat C\longrightarrow \t C$ be a semistable reduction of
 the pullback fibers $ \R_\pi$ of the branched fibers $\B_\pi$. Then we know
that $\pi\circ \hat\pi$ is also
the semistable reduction of $\B_\pi$.
By Lemma~2.4 and the basic equalities we can obtain the equalities in
this theorem.
\QED
\enddemo

\vskip0.8cm
\noindent
{\bf 3 \ On the invariants of non-semistable fibers}
\vskip0.6cm
\noindent
{\it 3.1 \ The computations of the invariants $c_1^2,c_2$ and $\chi$}
\vskip0.6cm
\noindent
In what follows, we shall consider the computation of the invariants
$c_1^2(F)$, $c_2(F)$ and $\chi_F$. By Noether equality, we only need
to compute $c_1^2$ and $c_2$.

First note that if we use embedded resolution to resolve the singularities
of $F$, then the number
    $$
   c_{-1}(F)={1\o d}\#\{ (-1)\text{-curves over $F'$ contracted by $\t\rho$}\}
    $$
is also independent of the stabilizing base change if $d$
is divided by $M_F$.

\proclaim{Theorem~3.1}
    $$ \eqalign{
        c^2_1(F)&   =4N_F+F^2_{\text{red}}
                       +\sum_{p\in F}\alpha_p-c_{-1}(F),  \cr
        c_2(F)  &   =2N_F+\mu_F-\sum_{p\in F}\beta_p  + c_{-1}(F).}
\eqno(15)
    $$
\endproclaim
\demo{Proof}
 The first formula has been proved in Proposition~2.2.
In order to prove the second formula, we consider the stabilizing base
change $\pi$ of $F$ whose degree is divided by $M_F$.
 By definition, if $\t F$ and $F_2$ are respectively the
pullback fibers of $F$ in $\t S$ and $S_2$, (note that $S_2$ is the
embedded resolution, not the minimal resolution),  then we have
    $$\eqalign{
       e_\pi&=e_F-{1\o d}e_{\t F}
             =e_F-{1\o d}e_{F_2}+c_{-1}(F). \cr}
    $$
Since $F_2$ is semistable, $e_{F_2}$ is the number of singular points
of $F_2$, which is exactly the number $d\sum_{p\in F}\beta_p$ (Lemma~1.4).
We have known that $e_F=2N_F+\mu_F$, hence the
second formula has been obtained.
\QED
\enddemo

\proclaim\nofrills{\it Remark. } \ {\rm From the formulas above and the
Noether formula,  we can see that
$\chi_F$ is independent of $c_{-1}(F)$, hence it can be computed directly
from embedded resolution. In fact, if we consider the canonical semistable
reduction of $F$, then we can prove that
    $$c_{-1}(F)= \sum_{q'\in F'}\beta_{q'},$$
where $F'$ is the embedded resolution of $F$, and $q'$ runs over the singular
points of $F'$ such that the $(-2)$-curves coming from $p'$ are contracted to
points
of the semistable model of $F$. }
\endproclaim

\proclaim\nofrills{\it Example. } \  {\rm Note that the discussion above
holds for elliptic fibrations. In this case, $K_\pi^2=0$ for all base changes,
so we have $c_1^2(F)=0$. By a direct computation we have
    $$c_2(F)=12\chi_F=\cases
              0, & \text{if $F$ is of type $_m$I$_b$}, \cr
              6, & \text{if $F$ is of type I$^*_b$ ($b>0$)}, \cr
              e_F, & \text{otherwise}.
                     \endcases
    $$}
\endproclaim

The result above shows the well-known fact that the semistable model of
an elliptic fiber is smooth except for type
$_m$I$_b$ ($b>0)$ and type I$^*_b$ ($b>0$).

\vskip0.8cm
\noindent
{\it 3.2 \ Proof of Theorem~C}
\vskip0.6cm
\noindent
\proclaim{Lemma~3.2}
     $$
     \sum_{p\in F}\alpha_p \leq 2p_a(F_{\text{red}})   ,            \eqno(16)
     $$
the equality holds iff $p_a(F_{\text{red}}) =0$, hence $F$ is a tree
of nonsingular rational curves.
\endproclaim
\demo{Proof} We shall use the notations of Sect.~1.1.
By (1) and (4), we have
     $$\eqalign{
           \alpha_p&=2\delta_p-\sum_{i=0}^{r-1} (m_i-2)-2(k_p-1)\cr
                   &\leq 2\delta_p-2(k_p-1).                 }       \eqno(17)
     $$
On the other hand, if the reduced part of $F$ consists of $l_F$ components
$\G_i$, $i=1,\cdots, l_F$, then
we have
     $$
           p_a(F_{\text{red}})=\sum_{i=1}^{l_F}p_a(\t \G_i)+
                   \sum_{p\in F}\delta_p -l_F+1.
                                                                    \eqno(18)
     $$
where $\t\G_i$ is the normalization of $\G_i$. Hence we only need to prove that
     $$
     \sum_{p\in F} (k_p-1) \geq l_F-1.                             \eqno(19)
     $$
But this inequality is an immediate consequence of the connectedness of $F$.
So (16) holds.

If the equality in (16) holds, then from (17) we know that $\alpha_p=0$ for
any $p\in F$, hence $p_a(F_{\text{red}})=0$. Then from (18) and (19),
we can see $F$ is a
tree of smooth rational curves.
\QED
\enddemo

\proclaim{Theorem~3.3}
    $$                  c_1^2(F)\leq 4g-4 .                           \eqno(20)
\noindent     $$
\endproclaim
\demo{Proof}
{}From Lemma~3.2 we have
    $$c_1^2(F)\leq 4N_F+F^2_{\text{red}}+2p_a(F_{\text{red}})-c_{-1}(F),
    $$
and the equality holds iff $p_a(F_{\text{red}})=0$.
Hence it is easy to prove that
    $$c_1^2(F)+c_{-1}(F)\leq 4g-3 ,$$
and the equality holds iff $F$ satisfies
    $$
          p_a(F_{\text{red}})=0, \hskip0.3cm F_{\text{red}}K=1  .  \eqno(21)
    $$
So it is enough to prove that for the fibers $F$ satisfying (21) we have
$c_{-1}(F)\geq 1$.

    Now we consider the canonical semistable reduction $\phi_F$ of $F$,
we know that the degree of $\phi_F$ is $M_F$.
We can see that the fiber $F$ satisfying (21) is a tree of a ($-3$)-curve
$\G$ and some $(-2)$-curves $E$. We note first that if $F$ contains a
($-2$)-curve $E$ such that  $F_{\text{red}}$ has only
one singular point $p$ on $E$,
then $p$ is an ordinary double point of type $(n, 2n)$, where $n$ is the
multiplicity of $E$ in $F$. Since the pullback fiber $\t F$ of $F$
is semistable, for any component $\G$ in $\t F$,
$-\G^2$ is the intersection number of $\G$ with the other components.
Thus we can see easily that the inverse image of $E$ in the minimal resolution
surface consists of
$n$ ($-1$)-curves, hence the exceptional curves of $p$ can be contracted to
a point. That is to say we contracted $n+[n, 2n]d-(n, 2n)=d/2$ \
$(-1)$-curves (Lemma~1.4).
Thus the contribution of $(E, p)$ to $c_{-1}(F)$ is
${1\o 2}$. On the other hand, we
know easily that there are at least two such $(-2)$-curves in $F$, hence
$c_{-1}(F)\geq 1$. This completes the proof.
\QED
\enddemo

\vskip0.8cm
\noindent
{\it 3.3 \ Proof of Theorem~B}
\vskip0.6cm
\noindent
\proclaim{Theorem~3.4}  For any singular fiber $F$, we have
       $$
                   c_1^2(F) \leq 2c_2(F),                      \eqno(22)
       $$
with equality iff $F=nF_{\text{red}}$ and $F$ has at worst ordinary
double points as its singularities.
\endproclaim
\demo{Proof} From (15), we have
    $$
2c_2(F)-c_1^2(F)=3c_{-1}(F)-F^2_{\text{red}}+\sum_p(2\mu_p-2\beta_p-\alpha_p).
    $$
Then by Lemma~1.3, $\mu_p-\beta_p\geq \alpha_p$, hence we have
    $$
       2c_2(F)-c_1^2(F) \geq -F^2_{\text{red}}+\sum_{p\in F}\alpha_p\geq 0.
    $$
If $c_1^2(F)=2c_2(F)$, then $F_{\text{red}}^2=0$, and $\alpha_p=0$ for
all $p\in F$. By the well-known Zariski's lemma [BPV, p.90],
we have $F=nF_{\text{red}}$.
Furthermore,  $\alpha_p=0$ implies that $p$ is an ordinary double point,
so $F_{\text{red}}$ has at worst nodes as its singularities.

The converse is obvious.
\QED
\enddemo

\proclaim{Proposition~3.5} If all of the multiple components of $F$ are
$(-2)$-curves,
 then
      $$
           c_1^2(F)\leq c_2(F).                              \eqno(23)
      $$
\endproclaim
\demo{Proof} The proof is similar to that of Theorem~3.4.
\QED
\enddemo

\vskip0.8cm
\noindent
{\bf 4 \ Applications}
\vskip0.6cm
\noindent
{\it 4.1 \ On the slopes of fibrations}
\vskip0.6cm
\noindent
{}From Theorem~3.4, Noether equality and the corollary to Theorem~A, we have

\proclaim{Theorem~4.1} For any stabilizing base change $\pi$, we have
      $$
                K_\pi^2  \leq 8 \chi_\pi.                      \eqno(24)
      $$
\endproclaim

As in the case of fibrations, we have the following definition of slopes.

\proclaim{Definition 4.2} {\rm If $F$ is a non-semistable fiber,
then $\chi_F\neq 0$, and so we can define the slope of $F$ as
         $$
            \lambda_F = {c_1^2(F)/ \chi_F}.
         $$
{}From Theorem~3.4, we know $0< \lambda_F\leq 8$.

If $\pi$ is a non-invariant base change, then we define the slope of $\pi$ as
         $$
            \lambda_\pi=K_\pi^2 /\chi_\pi .
         $$}
\endproclaim

Note that a non-trivial stabilizing base change $\pi$ satisfies $\chi_\pi>0$,
so
Theorem~4.1 says that its slope  $\lambda_\pi \leq 8$.

We have known in the introduction that for a stabilizing base change $\pi$,
       $$
      (\lambda_{\t
f}-\lambda_f)(\chi_f-\chi_{\B_\pi})=(\lambda_f-\lambda_\pi)\chi_{\B_\pi}.
                                                                   \eqno(25)
       $$

\proclaim{Corollary~4.3} If $f:S\longrightarrow C$ is a non-semistable
fibration
with $\lambda_f>8$, then through any non-trivial stabilizing base change,
we have
       $$
              \lambda_{\t f}>\lambda_{f}.                       \eqno(26)
       $$
\endproclaim

In what follows, we shall consider a set of fibrations $\Sigma$
which is invariant under base changes, i.e., if $f\in \Sigma$, then
$\t f\in \Sigma$.

\proclaim{Corollary~4.4} Let $f$ (resp. $f'$) be a fibration in $\Sigma$
with maximal (resp. minimal) slope.

     {\rm 1)} For any non-semistable fiber $F$ of $f$
(resp. $f'$), we have
        $$
             \lambda_F\geq \lambda_f, \hskip1cm
              (\text{resp. }  \lambda_F\leq \lambda_{f'}).
        $$

     {\rm 2)} If $\lambda_f > 8$, then $f$ is semistable.

     {\rm 3)} If $\lambda_f>6, $ then any non-semistable fiber of $f$ has at
least one multiple component which is not a $(-2)$-curve.
\endproclaim
\demo{Proof} Considering the canonical stabilizing base change of $F$
and using (25), we can prove 1).  2) and 3) are immediate consequences
of (22)--(25) and the assumption.
\QED
\enddemo

\proclaim\nofrills{\it Remark. } \  {\rm This corollary can be used to classify
singular fibers of a fibration with minimal slope in the sense above. For
example,

    I)  Xiao [Xi1, Xi4] has proved that for any relatively minimal
fibration $f$ of genus $g$,
       $$\lambda_f\geq 4-4/g.$$
Furthermore, if $f$ is a hyperelliptic fibration, then
        $$\lambda_f\leq 12-{4g+2\o [g^2/2]}.$$

    II) If $f$ is non-hyperelliptic, then the lower bounds $\lambda_g$ of the
slope are $\lambda_3=3$, $\lambda_4=24/7$,  $\lambda_5=40/11$.
(cf. [Ch, Ho, Ko, Re]).
}\endproclaim

If we consider fibrations over $\Bbb P^1$, and we only consider
base changes with two ramification points, then the above results
can also be used.

\vskip0.8cm
\noindent
{\it 4.2 \ Canonical class inequality for general fibrations}
\vskip0.6cm
\noindent
First we recall Miyaoka's inequality and refer to [Hi] for the details.

\proclaim{Lemma~4.5} {\rm [Mi]} If $S$ is a smooth minimal surface of general
type, and $E_1$, $\cdots$, $E_n$ are disjoint $ADE$ curves on $S$, then
we have
     $$
       \sum_{i=1}^nm(E_i)\leq 3c_2(S)-c_1^2(S),
     $$
where $m(E)$ is defined as follows,
     $$
       \eqalign{ m(A_r)&=3(r+1)-{3\over r+1} ; \hskip0.3cm
                 m(D_r)=3(r+1)-{3\over 4(r-2)}, \hskip0.2cm \text{ for }
r\geq4;                    \cr
                 m(E_6)&=21-{1\over8};     \hskip0.3cm
                 m(E_7)=24-{1\over16};               \hskip0.3cm
                 m(E_8)=27-{1\over40}.                     \cr}
    $$
\endproclaim

The condition ``of general type'' can be replaced by some other conditions.
(cf. [Mi]).

\proclaim{Theorem~4.6} If $f$ is a fibration of genus $g>1$ over a curve $C$
of genus $b$, then
    $$
     K_{S/C}^2\leq 3\sum_{y\in C} \delta_y^\#+(2g-2)\max(2b-2, 0),
\eqno(27)
    $$
where $\delta_y^\#=e_{F_y}-{1\over 3}\sum_{E\subset F_y}m(E) \leq 4g-3$, and
the
sum is taken over all of the disjoint ADE curves $E$ in $F_y$.
\endproclaim
\demo{Proof} Inequality (27) is an immediate consequence of
Miyaoka-Yau inequality (See Vojta's proof, [Vo, Theorem~2.1]).
So we only need to prove that $\delta_y^\#\leq 4g-3$.

  We denote by $l_D$ the number of components of a curve $D$, and by $\t D$
the normalization of $D$. Let
$F_{\text{red}}=D+\sum_{E\subset F}E$ be the reduced part of a fiber $F$.
 Then
  $$
       \eqalign{
      \chi_{\text{top}}(F)&=\c_t(D)+\sum_{E\subset F}\left(\c_t(E)
           -\#(D\cap E)\right) \cr
                          &=\c_t(\tilde D)-\sum_{p\in D}(k_p(D)-1)
                  +\sum_{E\subset F}(l_E+1)-\#(D\cap\sum_{E\subset F}E)
                                                           \cr
                          &\leq 2l_D+\sum_{E\subset F}(l_E+1)
                 -\sum_{p\in D}(k_p(D)-1)-\sum_{E\subset F}\#(D\cap E). }
   $$
{}From the definition of $m(E)$, we know
        $$m(E)\geq 3(l_E+1)-1, $$
hence
   $$ \eqalign
      {  \delta^\#_{F}&=2g-2+\c_t(F)-\sum_{E\subset F}m(E)    \cr
                     &\leq 2g-2 +2l_D-\sum_{p\in D}(k_p(D)-1)-\#(D\cap
            \sum_{E\subset
                       F}E)+\sum_{E\subset F}1.\cr}
   $$
Since $l_D\leq KD\leq 2g-2$, it is enough to prove that
   $$ \sum_{p\in D}(k_p(D)-1)+\sum_{E\subset F}(\#(D\cap E)-1)\geq l_D-1.
                                                              \eqno(28)
   $$
Indeed, if $D$ is connected, then we have
   $$
     \sum_{p\in D}(k_p(D)-1)\geq l_D-1,
   $$
hence (28) holds. If $D$ has $r$ connected components $D_1$, $\cdots, $$D_r$,
then
   $$\sum_{p\in D}(k_p(D)-1)\geq \sum_{i=1}^r(l_{D_i}-1)=l_D-r.$$
On the other hand, from the connectedness of $F$, we can prove easily that
   $$\sum_{E\subset F}(\#(D\cap E)-1)\geq r-1.$$
Hence (28) holds too.
\QED
\enddemo

In [Vo], Vojta shows that if $f$ is a non-trivial semistable fibration with $s$
singular fibers, then

    $$
         K^2_{S/C}\leq (2g-2)(2b-2+s).
    $$
By using semistable reductions, we can obtain a similar inequality for general
 fibrations.

\proclaim{Theorem~4.7} If $f$ is a non-trivial fibration of genus $g\geq 2$
 with $s$ singular fibers, then
    $$
    K_{S/C}^2\leq (2g-2)(2b-2+3s),                                  \eqno(29)
    $$
and the equality holds only if $f$ is smooth.
\endproclaim
\demo{Proof} We denote by $F_1,\cdots,F_s$ the $s$ singular fibers of $f$,
and assume that $s>0$. Note that $f$ is non-trivial, so $b=0$ implies
that $s\geq 2$.
In fact, if $b=0$ and $s=2$, then  $f$ is isotrivial (cf. [Be1]), thus
 $I_K(f)=0$. In this case, the inequality (20) holds strictly, so does (29).
Hence we can assume that $s\geq 3$ if $b=0$.

First we claim that there exist some semistable reductions $\pi:\wt
C\longrightarrow C$ of $f$ satisfying the following two conditions:

    {\rm(i)} $\pi$ is ramified uniformly over the $s$ critical points of $f$,
and the
ramification index of $\pi$ at any ramified point is exactly $e$.

    {\rm (ii)} $e$ is divided by $M_{F_i}$ for all $i$, and it can be
arbitrarily large.

In fact, a base change satisfying the above two conditions must be a semistable
reduction. If $b=g(C)>0$, then the existence follows from Kodaira-Parshin's
construction.
If $b=0$, then $s\geq 3$ by assumption.
Hence we can construct a base change totally ramified over the $s$ points.
Then the existence is reduced to the case $b>0$.

We let $\wt f: \wt S\longrightarrow \wt C$ be the semistable model of $f$ under
the base change $\pi$ constructed above, $\wt s$ the number of singular fibers
of $\wt f$,
and $\wt b$ the genus of $\wt C$. Denote by $d$ the degree of $\pi$.
Then we know that $\wt s\leq ds/e$, and
    $$
     2\wt b-2=d(2b-2)+d\left(1-{1\over e}\right)s.
    $$
{}From Theorem~A, we obtain
    $$
       K_f^2={1\over d}K^2_{\wt f}+\sum_{i=1}^sc_1^2(F_i),
    $$
hence we have
    $$\eqalign{
     K_f^2-(2g-2)(2b-2+3s)
    =&{1\over d}\left(K^2_{\wt f}-(2g-2)(2\wt b-2+\wt s)\right) \cr
     +{2g-2\over d}\left(\wt s-{d\over e}s\right)
      +&\sum_{i=1}^s\left(c_1^2(F_i)-(4g-4)\right).}
    $$
Consequently, combining Theorem~3.3 and Vojta's canonical class inequality for
semistable fibrations, we can obtain immediately (29).
Since the equality in Vojta's inequality implies the smoothness
of the fibration (cf. [Ta2, Lemma~3.1]), we can see easily that if
the equality in (29) holds, then $f$ is smooth.
\QED
\enddemo

\vskip0.8cm
\noindent
{\it 4.3 \ On Horikawa number of a non-semistable fiber of genus 3}
\vskip0.6cm
\noindent
Let $f:S\longrightarrow C$ be a relatively minimal non-hyperelliptic fibration
of genus 3,
$F$ a fiber of $f$, and $p=f(F)$. Then the Horikawa number of $F$ is
defined as (cf. [Re])
    $$
     H_F=\text{length coker}\left(S^2f_*\omega_{S/C}\hookrightarrow
      f_*(\omega_{S/C}^{\otimes 2})\right)_p.
    $$
The global invariants of $f$ depend on this number. In fact, Reid [Re] shows
that
     $$K_f^2-3\chi_f=\sum_{F}H_F .                               \eqno(30)
      $$
In general, it is quite difficult to compute $H_F$. The aim here is to
try to reduce the computation for a non-semistable fiber to the
computation  for its semistable models,
by using semistable reductions.

\proclaim{Theorem~4.8} Let $\t F$ be the semistable model of $F$ under a
totally ramified stabilizing base change of degree $d$. Then
    $$
           {1\o d}H_{\t F}=H_F+{1\o 4}(c_2(F)-3c_1^2(F)).    \eqno(31)
    $$
\endproclaim
\demo{Proof} We can assume that the branch locus of the base change consists of
$F$ and some generic smooth fibers, hence
    $$\eqalign{c_1^2(F)-3\chi_F&=K_\pi^2-3\chi_\pi\cr
                  &=K_f^2-3\chi_f-{1\o d}(K_{\t f}-3\chi_{\t f})\cr
            &=H_F-{1\o d}H_{\t F}.}
    $$
By using Noether formula, we can obtain (31).
\QED
\enddemo

\proclaim\nofrills{\it Examples}. \ {\rm If $F$ is a genus 2 curve with an
ordinary cusp, then we can take $d=6$. Then the semistable model
$\t F$ consists of a nonsingular elliptic curve  $E$
and a nonsingular curve $C$ of genus 2, with $EC=1$.
Since $c_1^2(F)={1\o 6}$ and $c_2(F)={11\o6}$,
 we have
    $$H_{\t F}=6H_F+2.$$

      If $F=2C$, $C$ is a smooth curve of genus 2, then we can take $d=2$,
hence $\t F$ is a nonsingular hyperelliptic curve of genus 3. We have
     $c_1^2(F)=4, \ c_2(F)=2$.
Hence
    $$H_{\t F}=2H_F-5.$$
So we can compute directly the Horikawa numbers of some special
singular fibers, e.g., if their semistable models are non-hyperelliptic curves
of genus 3.}
\endproclaim

\vskip0.8cm
\noindent
{\bf 5 \  The proof of Lemma~1.5}
\vskip0.6cm
\noindent
In this section, we shall use freely the notations of Sect.~1.2.
Note first that Lemma~1.5 is a special case of the following theorem.

\proclaim{Theorem~5.1} For the embedded resolution given in Sect.~1.2,
we have
    $$
    -K^2_{\phi}=d\sum_{i=0}^{r-1}\left(m_i-2+{1\over d}\left((m_i^*, d)-m_i(d)
\right) \right)^2-\sum_{q\in B_r}K_q^2 ,
    $$
where $m_i, m_i^*, m_i(d)$ are respectively the multiplicities of
$B_{i, \text{red}}$, $B_i$, $B_i(d)$ at $p_i$, and
$B_i(d)=\sum_{\Gamma}(d, n_\G)\G$ if $B_i=\sum_{\G} n_\G \G$.
\endproclaim
\demo{Proof}
Since we only need to find $K_\phi$, without loss of generality, we
may assume that $U_0$ is a compact smooth surface, and
the reduced curve of $B=B_0=\sum_{\G}n_\G \G$ has only one singular point $p$,
(otherwise we can resolve in advance the other singularities of $B$ by using
embedded
resolution).  So we have a formula similar to (11):
   $$
    K_M=\phi^*\pi_0^*\left( K_{U_0}+\sum_{\G\subset B}
       \left(1-{(d, n_\G)\o d}\right)\G\right)+K_\phi,
    $$
i.e.,
    $$
          K_M=\eta^*\pi^*\sigma^*\left(K_{U_0}+B_{0, \text{red}}
          -{1\o d}B_0(d)\right)+K_\phi.                              \eqno(32)
    $$
On the other hand, $\pi_r$ is determined by $B_r=\sigma^*(B)$, hence we have
also
     $$
      K_M=\eta^*\pi_r^*\left(K_{U_r}+B_{r, \text{red}}-{1\o d}B_{r}(d)\right)
          +K_{\eta}   .                                      \eqno(33)
     $$
{}From Definition 1.1 it is easy to prove that
   $$\eqalign{
           K_{U_i}+B_{i, \text{red}}-{1\o d}B_i(d)=&\sigma_i^*\left(
      K_{U_{i-1}}+B_{i-1, \text{red}}-{1\o d}B_{i-1}(d)\right)
                   \cr &   -\left(m_{i-1}-2+{1\o d}((m^*_{i-1}, d)-m_{i-1}(d))
                         \right)E_i.}                         \eqno(34)
    $$
If we denote by $\Cal E_i$ the total inverse image of $E_i$ in $U_r$, then
from (32)--(34), we have
   $$
      K_\phi=-\eta^*\pi^*_r\left(\sum_{i=1}^r\left(m_{i-1}-2+{1\o d}(
      ( m_{i-1}^*, d)-m_{i-1}(d))\right)\Cal E_i\right) +K_\eta.
   $$
Hence we obtain the desired equality.
\QED
\enddemo

\proclaim\nofrills{\it Remark}. \ {\rm In order to resolve the singularities of
 $V_0$,
we can also use the {\it $d$-resolution of $(B,p)$}, which is defined as in
Definition~1.1 by replacing condition (3) with

($3^\prime$)   If $E_i=\sigma_i^{-1}(p_{i-1})$ is the exceptional curve, then
we have
        $$
             B_i=\sigma_i^*(B_{i-1})-d\left[{m_{i-1}^*\over d}\right]E_i .
        $$
Then we have also the formula in Theorem~5.1.}
\endproclaim

Finally, we consider the computation of $K^2_\eta$ by using Jung-Hirzebruch's
resolution.

Let $V$ be the normalization of the local surface defined by $z^d=x^ay^b$, and
$\eta:M\longrightarrow V$ the minimal resolution of the singularities. Then
the resolution can be constructed as follows (cf. [BPV,p.83]).

    I. \ If $(d, a)=(d, b)=(a, b)=1$, then  we let $q, q'$ be two integers with
$0<q, q'<d$, $aq+b\equiv 0\pmod d$, $qq'\equiv 1\pmod d$. If
      $${d\over q}=e_1-\cfrac1\\e_2-\cfrac1\\
                \ddots -\cfrac1\\ e_r\endcfrac\ ,
      $$
then the exceptional set of the resolution is a chain of $r$ rational curves
      $$
   {\overset{-e_1}\to\bullet}\hskip-.2cm{\hskip0.8cm\over\hskip1cm}\hskip-0.2cm
      {
\overset{-e_2}\to\bullet}\hskip-.2cm{\hskip0.8cm\over\hskip1cm}\hskip-0.2cm
{\overset{}\to\bullet}\hskip-.2cm{\hskip0.8cm\over\hskip1cm}\hskip0.3cm\cdots
   \hskip0.3cm    {
\overset{}\to\bullet}\hskip-.2cm{\hskip1cm\over\hskip1cm}\hskip-0.2cm
       \overset{-e_r}\to\bullet
     $$
and we have
    $$
    -K_\eta^2=\sum_{i=1}^r(e_i-2)+{q+q'+2\over d}-2.
    $$

   II. \ If $(d, a, b)=1$, then the singularity of $V$ is isomorphic
to the normalization of $z^{d'}=x^{a'}y^{b'}$, where
$a=a'(d, a)$, $b=b'(d, b)$ and $d=d'(d, a)(d, b)$, hence the resolution and the
computation are reduced to (I).

   III. \ If $d_0=(d, a, b)>1$, then we have
      $$
    z^d-x^ay^b=\prod_{i=1}^{d_0}\left(z^{d/d_0}-x^{a/d_0}y^{b/d_0}\exp\left
(2\pi i\sqrt{-1}
/d_0\right)\right).
      $$
Hence the singularity decomposes into $d_0$ singularities of type II.

In particular, if $d$ is divided by $a$ and $b$, by (III) and (II), we can see
that $V$ consists of $d_0=(a,b)$ singular points of type $A_n$,
where
    $$
         n={d_0d\over ab}-1.
    $$
This is what we expected in Lemma~1.4.
\vskip0.4cm
\proclaim\nofrills{\eightit Acknowledgement. } \
{\eightrm  The author would like to thank Prof. Xiao Gang for encouraging me to
find the
best possible inequalities between the invariants of base changes. The author
thanks
also the referee for the valuable suggestion for correction of the original
version.}
\endproclaim

\enddocument